\documentclass[11pt]{article}

\usepackage{amsmath}
\usepackage{amscd}
\usepackage{amsfonts}
\usepackage{epsfig}
\usepackage{theorem}
\usepackage{authblk}
\usepackage[titletoc]{appendix}
\usepackage[margin=1.0in]{geometry}

\def \C {\mathbb C}

\def \R {\mathbb R}
\def \N {\mathbb N}
\def \Z{\mathbb Z}

\title{A Generalized Lerche-Newberger Formula}

\author{Parker Kuklinski, Michael Warnock, David Hague}

\begin{document}

%%%%%%%%%%%%%%%%%%%%%%%%%%%%%%%%%%%%%%%%%%%%%%%%%%%%%%%%%%%%%%%%%%%%%%%%%%
%setup commands for the bu thesis style file

\maketitle

\begin{abstract}
The Lerche-Newberger formula simplifies harmonic sums of Bessel functions and has seen application in plasma physics and frequency modulated quantum systems. In this paper, we rigorously prove the formula and extend the classical result to a family of multi-dimensional extensions of the single variable Bessel functions called generalized Bessel functions. Since prevailing definitions of these functions do not accommodate arbitrary complex order, we use an auxiliary family of functions called generalized Anger functions and show that the single-variable result holds in multiple dimensions for a certain selection of parameters. We conclude by applying these results to physical systems.
\end{abstract}

\section{Introduction}

Sums involving Bessel functions are pervasive in physics \cite{lerche08} \cite{szapiel78} \cite{watson22}. A Jacobi-Anger expansion, for instance, allows for a Fourier representation of a sinusoidal frequency modulated signal as a sum of harmonic components weighted by Bessel functions of varying order \cite{abramowitz64}. Other sums like the Kapteyn series arise in astrophysics applications \cite{lerche07}, while Neumann and Schl\"{o}milch-type series also also appear in the literature \cite{kravchenko17} \cite{linton06}. The Lerche-Newberger formula \cite{newberger82} gives a closed form representation to an infinite sum of Bessel function products with a transposed harmonic term (i.e. $1/(n-a)$). These sums first arose in 3-dimensional plasma systems with an oscillating ambient magnetic field from the associated $3\times 3$ plasma wave dispersion relation \cite{buneman59} \cite{lerche08a}. In the present context we are interested in applying the Lerch-Newberger formula to quantum systems under frequency modulation \cite{silveri17}. Berns et. al. \cite{berns06} investigate an intermediate ''quasiclassical" regime in a driven two level quantum system, and the transition rate of the qubit is modeled by a modified Lerche-Newberger summation.

In this paper we aim to adapt the methods of Kuklinski and Hague \cite{kuklinski19} to extend the traditional Lerche-Newberger formula to a class of multi-dimensional Bessel functions called generalized Bessel functions \cite{dattoli96}. Generalized Bessel functions introduce higher order harmonics in the modulation function of the associated integral. In a signal processing context, the generalized Bessel functions are the Fourier coefficients of a multi-tone sinusoidal frequency modulated (MTSFM) signal; these signals are useful for radar and sonar applications due to their constant amplitude, spectral efficiency, and tunability \cite{Hague_AES}. The parameters of these MTSFM are fed through to the arguments of the GBFs in the Fourier transformed space. The generalized Bessel functions have found applications in laser physics \cite{reiss62}, crystallography \cite{paciorek94}, and astrophysics \cite{korsch06}. In connection to frequency modulated quantum systems, a version of the Lerche-Newberger sum found in Berns et. al. with Bessel functions replaced by GBFs arises in a more general oscillatory system.

A na\"{i}ve extension of the Lerche-Newberger formula by a total substitution of one-dimensional Bessel functions with their multi-dimensional counterparts fails since the generalized Bessel functions are not well defined for fractional order. To remedy this, we instead work with generalized Anger functions which extend the usual integral representation of integer-order generalized Bessel functions to arbitrary order. These Anger functions, although agreeing with Bessel functions on integer order, do not satisfy the same collection of identities and thus certain unresolvable terms emerge in the corresponding Lerche-Newberger extension.

The rest of this paper is organized as follows: In section 2 we rigorously derive the one-dimensional Lerche-Newberger formula over the usual domain of parameters. In section 3 these results are extended by replacing the one-dimensional Bessel functions with generalized Anger functions. Section 4 treats an application to a multiply frequency modulated two-level quantum system. An appendix handles some details of the proof of the one-dimensional Lerche-Newberger formula.

\section{Derivation of Lerche-Newberger formula}

In this section we will derive the following sum presented by Newberger \cite{bakker84} \cite{newberger82}:
\begin{equation}
\sum _{n=-\infty}^\infty\frac{(-1)^nJ_{\alpha -\gamma n}(z)J_{\beta +\gamma n}(z)}{n+\mu}=\frac{\pi}{\sin \mu\pi}J_{\alpha +\gamma\mu}(z)J_{\beta -\gamma\mu}(z)
\end{equation}
Here, we restrict $\text{Re}(\alpha +\beta )>-1$, $\mu\in\C\backslash\Z$, and $\gamma\in (0,1]$. To proceed, we follow an argument from Lerche \cite{lerche74}; recall from Watson \cite{watson22} that for $\text{Re}(\alpha +\beta)>-1$ we have
\begin{equation}J_\alpha (z)J_\beta (z)=\frac{2}{\pi}\int _{0}^{\pi /2} J_{\alpha +\beta}(2z\cos\theta )\cos (\alpha -\beta )\theta d\theta .
\end{equation}
Substituting this into the summation on the left side of (1), we have
\begin{equation}
\sum _{n=-\infty}^\infty\frac{(-1)^nJ_{\alpha -\gamma n}(z)J_{\beta +\gamma n}(z)}{n+\mu}=\frac{2}{\pi}\sum _{n=-\infty}^\infty\left[\int _0^{\pi /2}\frac{(-1)^n}{n+\mu}J_{\alpha +\beta}(2z\cos\theta )\cos{(\alpha -\beta-2\gamma n)\theta}d\theta\right] .
\end{equation}
Ultimately we would like to interchange summation and integration, and we do this with dominated convergence theorem \cite{royden88}. Since the doubly-infinite summation is a composition of two separate limiting operations, one over the negative indices and one over the positive indices, we fmust split the double summation into two one-sided summations. one we label $S_+$ with $n\in\{ 1,2,...\}$ and the other we call $S_-$ with $n\in\{ -1,-2,...\}$ such that
\begin{equation}
\sum _{n=-\infty}^\infty\frac{(-1)^nJ_{\alpha -\gamma n}(z)J_{\beta +\gamma n}(z)}{n+\mu}=\frac{2}{\pi}\left( S_++S_-+\frac{1}{\mu}\int _0^{\pi /2}J_{\alpha +\beta}(2z\cos{\theta})\cos{(\alpha -\beta )\theta}d\theta\right) .
\end{equation}

To apply dominated convergence theorem to the summations $S_\pm$, we first state a result for two simpler sequence of functions. Let $f_n(\theta )$ and $h_n(\theta )$ be the sequences of partial sums
\begin{equation}
f_n(\theta )=\sum _{k=1}^n\frac{(-1)^k}{k+\mu}\cos (k\theta ),\hspace{1cm}g_n(\theta )=\sum _{k=1}^n\frac{(-1)^k}{k+\mu}\sin (k\theta )
\end{equation}
for $\theta$ in some finite interval in $\R$. These partial sums are related to Lerch zeta functions \cite{apostol51}. We are searching for functions $F(\theta )$ and $G(\theta )$ independent of $n$ such that $|f_n(\theta )|<F(\theta )$ and $|g_n(\theta )|<G(\theta )$ for all $n$ and $F(\theta )$ and $G(\theta )$ are integrable. We can use a triangle inequality to bound $|f_n(\theta )|$ in the following way:
\begin{align}
|f_n(\theta )| &\le \left\lvert\sum _{k=1}^n\frac{(-1)^k}{k+\mu}\cos (k\theta ) -\sum _{k=1}^n\frac{(-1)^k}{k}\cos (k\theta )\right\rvert +\left\lvert\sum _{k=1}^n\frac{(-1)^k}{k}\cos (k\theta )\right\rvert \nonumber \\
	&= |\mu |\left\lvert\sum _{k=1}^n\frac{(-1)^k}{k(k+\mu )}\cos (k\theta )\right\rvert +\left\lvert\sum _{k=1}^n\frac{(-1)^k}{k}\cos (k\theta )\right\rvert
\end{align}
If $\mu$ isn't a negative integer, then by a $p$-series test the first sum on the right hand side pf (6) can be bounded by a finite quantity. The second sum on the left is a partial sum of a logarithm, specifically on the open interval $(-\pi ,\pi )$ we have the pointwise convergence
\begin{equation}
\sum _{k=1}^\infty\frac{(-1)^k}{k}\cos (k\theta )=-\log\left(\cos\frac{\theta}{2}\right)-\log{2}
\end{equation}

We prove that the partial sum term in (6) is bounded above by its pointwise limit plus a constant, and is therefore bounded by an integrable function letting us use dominated convergence theorem. We unfortunately cannot apply standard results in Gibbs phenomenon \cite{wilbraham48} as these results apply only to discontinuous functions with finite jumps; the one we have here is an infinite discontinuity. Nevertheless, we can exploit the specific nature of our function to arrive at an appropriate bound. Let $s_n(\theta )$ be the $n^\text{th}$ partial sum such that $s_n(\theta )+\log\left(\cos\theta /2\right)+\log{2}$ has local extrema in this interval at points $\theta _j=2\pi j/(2n+1)$. We prove in the appendix that there exists some $M\in\N$ such that for all $n>M$ the largest extrema of this difference occurs at $\theta _n=\pi-\pi /(2n+1)$ (that this holds only for $n$ sufficiently large does not impede the conditions of the DCT since for all $n<M$ we can bound (6) by a constant). Therefore by plugging $\theta _n$ into the difference of these terms we have
\begin{equation}
s_n(\theta _n)+\log\left(\cos\frac{\theta _n}{2}\right)+\log{2}=\left[\sum _{k=1}^n\frac{1}{k}\cos\left( \frac{k\pi}{2n+1}\right)\right] +\log\left(\sin\frac{\pi}{4n+2}\right) +\log 2
\end{equation}
We show that as $n$ gets large, both the partial sum term and the logarithm term grow at $O(\log{n})$ with opposite leading coefficients, thus cancelling and leaving an O(1) term. To handle the partial sum term, use a trigonometric identity to extract a harmonic sum:
\begin{equation}
\sum _{k=1}^n\frac{1}{k}\cos\left( \frac{k\pi}{2n+1}\right)=\sum _{k=1}^n\frac{1}{k} -\sum _{k=1}^n\frac{2}{k}\sin ^2\left(\frac{\pi k}{4n+2}\right)
\end{equation}
By applying a standard result on harmonic series and using the identity $\sin x<x$ for all $x>0$, we see that this partial sum on the left in (9) is equal to $\log{n}+O(1)$. Next, we can conduct an asymptotic expansion of the logarithm term in (8). Note that for small $x$ we have $\sin{x}=x(1+O(x^2))$ such that the logarithm pulls out this factor of $x$ as well as other multiplicative factors and therefore:
\begin{equation}
\log\left(\sin\frac{\pi}{4n+2}\right) =-\log{n}+O(1)
\end{equation}

The two asymptotic expansions we've conducted imply that the quantity in (8) converges to a limit, or the partial sum in (6) overshoots its limit by an asymptotically finite quantity, and therefore there exists some finite $M$ such that
\begin{equation}
\left\lvert\sum _{k=1}^n\frac{(-1)^k}{k}\cos (k\theta )\right\rvert\le\frac{1}{2}\left\lvert\log\left( 1+\cos\theta\right) +\log{2}\right\rvert +M
\end{equation}
for all $n\in\N$ (Numerical simulations appear to show $M=1/2$ is a sufficient choice). This function on the right is Lebesgue integrable on any interval in $\R$ even in those that include singularities at $\theta =(2N+1)\pi$. Therefore, we can apply dominated convergence theorem to this sequence of partial sums $f_n(\theta )$. Applying dominated convergence theorem to $\{ g_n(\theta )\}$ is similar, but instead of bounding the alternating cosine summation by a diverging yet integrable function, we can bound the corresponding alternating sine function by a constant since this is the Fourier transform of the sawtooth wave.

To extend this to a partial sum of interest in $S_+$, we need to bound $J_{\alpha +\beta}(2z\cos\theta )[\cos (\alpha -\beta )\theta f_n(2\gamma\theta )+\sin (\alpha -\beta )\theta g_n(2\gamma\theta )]$ by a single integrable function for all $n$. But since Bessel functions and sine and cosine functions are entire on $\C$, we can bound these and $g_n(2\gamma\theta )$ by constants, and we can bound $f_n(2\gamma\theta )$ by the scaled version of the dominating function on the right hand side of (11). Thus, the interchange of summation and integration is justified. These arguments are directly applicable to $S_-$, however we also need to restrict $\mu$ from being a positive integer, and since a $1/\mu$ term appears on the right hand side of (4), we insist that $\mu\in\C\backslash\Z$.

Because dominated convergence theorem holds separately for $S_\pm$, there is no issue with now interchanging the entire double summation with integration in (3). We use the following identities to reduce this operation to a closed form
\begin{equation}
\sum _{n=-\infty}^\infty\frac{(-1)^n}{n+\mu}\cos (n\theta )=\frac{\pi\cos{\mu\theta}}{\sin{\pi\mu}},\hspace{1cm}\sum _{n=-\infty}^\infty\frac{(-1)^n}{n+\mu}\sin (n\theta )=-\frac{\pi\sin{\mu\theta}}{\sin{\pi\mu}}
\end{equation}
where $\theta\in [-\pi ,\pi ]$. In other words, the left hand side summations are the Fourier series for non-smooth $2\pi$-periodic functions. To apply these formulas to (4), we require that $\gamma\in (0,1]$ such that $2\gamma\theta\in [-\pi ,\pi]$ for the region of integration $\theta\in [0,\pi /2]$. Therefore we write the interchange as
\begin{equation}
\frac{2}{\pi}\sum _{n=-\infty}^\infty\left[\int _0^{\pi /2}\frac{(-1)^n}{n+\mu}J_{\alpha +\beta}(2z\cos\theta )\cos{(\alpha -\beta -2\gamma n)\theta}d\theta\right] =\frac{2}{\sin{\pi\mu}}\int _0^{\pi /2}J_{\alpha +\beta}(2z\cos\theta )\cos{(\alpha -\beta +2\gamma\mu )\theta}d\theta
\end{equation}
Since the right hand side is in the integral form of the Bessel product formula in (2), we can reverse the equation to arrive at the final result in (1). 

\section{Extension to generalized Anger functions}

In this section we create a generalized version of the Lerche-Newberger formula in (1) extended to generalized Anger functions. Traditionally, the two-dimensional generalized Bessel function is defined as the convolutional summation:
\begin{equation}
J_n(x,y)=\sum _{k=-\infty}^\infty J_{n-2k}(x)J_k(y)
\end{equation}
We usually call this function the index $(1,2)$ generalized Bessel function. If we naively attempt to apply this definition to the Lerche-Newberger formula in (1), we will run into problems because this form of the GBF is not defined for $n\in\C\backslash\Z$ as the summation in (14) will not converge. We could try to define the GBF in a completely analogous way as the Bessel functions by first defining a partial differential equation and defining the GBF as solutions to this PDF, but as seen in Kuklinski and Hague \cite{kuklinski19}, the index $(1,2)$ GBF satisfies two independent second-order linear PDEs, neither of which have clear extensions to GBFs of different index or higher order.

Regardless of these issues, however, if we restrict a putative Lerche-Newberger formula to GBFs of integer order, we will have a well-defined convergent sum. To do this, we introduce generalized Anger functions which are well defined for all orders but agree with the GBF for integer order. For two finite coordinates ${\bf x}=(x_1,...,x_m)$ and ${\bf y}=(y_1,...,y_m)$ in $\C^m$, we define the generalized Anger function of order $\alpha$ as:
\begin{equation}
A_\alpha ({\bf x},{\bf y})=\frac{1}{2\pi}\int _{-\pi}^{\pi}\exp\left[ i\left(\alpha\theta -\sum _{k=1}^m\left( x_k\sin{k\theta}+y_k\cos{k\theta}\right)\right)\right] d\theta
\end{equation}
These functions have favorable decay properties in the sense that a stationary phase approximation shows that for fixed ${\bf x},{\bf y}\in\C ^m$, $A_\alpha ({\bf x},{\bf y})=O(\alpha ^{-k})$ for all $k>0$ as $|\alpha |\rightarrow\infty$.

With these functions, we prove a new version of the Lerche-Newberger formula:
\begin{equation}
\sum _{n=-\infty}^\infty\frac{(-1)^nA_{\alpha -\gamma n}({\bf x},{\bf y})A_{\beta +\gamma n}({\bf x},{\bf y})}{n+\mu}=\frac{\pi}{\sin{\mu\pi}}A_{\alpha +\gamma\mu}({\bf x},{\bf y})A_{\beta -\gamma\mu}({\bf x},{\bf y})
\end{equation}
This equation holds for all $\mu\in\C\backslash\Z$, and all $\alpha ,\beta\in\C$ if $\gamma\in [0,1/2]$. We will prove this by breaking each of the generalized Anger functions into two integrals and passing the sum through using convergence theorems. Let us represent the first generalized Anger function as an integral:
\begin{align}
\sum _{n=-\infty}^\infty\frac{(-1)^nA_{\alpha -\gamma n}({\bf x},{\bf y})A_{\beta +\gamma n}({\bf x},{\bf y})}{n+\mu}=\frac{1}{2\pi}\sum _{n=-\infty}^\infty\int_{-\pi}^\pi &\frac{(-1)^ne^{-i\gamma n\theta}A_{\beta +\gamma n}({\bf x},{\bf y})}{n+\mu}\times ... \nonumber \\
	&...\times\exp\left[ i\left(\alpha\theta -\sum _{k=1}^m\left( x_k\sin{k\theta}+y_k\cos{k\theta}\right)\right)\right] d\theta
\end{align}
From the decay properties of the generalized Anger functions, we can use Fubini's theorem to interchange the summation and integral. Upon pulling the summation inside the first integral, we then expand $A_{\beta +\gamma n}({\bf x},{\bf y})$ into its integral representation and focus on the following quantity:
\begin{equation}
\sum _{n=-\infty}^\infty\frac{(-1)^ne^{-i\gamma n\theta}A_{\beta +\gamma n}({\bf x},{\bf y})}{n+\mu}=\frac{1}{2\pi}\sum _{n=-\infty}^\infty\int _{-\pi}^\pi\frac{(-1)^ne^{i\gamma n(\phi -\theta )}}{n+\mu}\exp\left[ i\left(\beta\phi -\sum _{k=1}^m\left( x_k\sin{k\phi}+y_k\cos{k\phi}\right)\right)\right] d\phi
\end{equation}
We can repeat the same argument from the previous section to interchange the summation and integral in (18).

However, when executing the complex exponential version of the sums described by (12), we must be careful since we do not necessarily have $\gamma(\phi -\theta )\in [-\pi ,\pi ]$ for general $\gamma$. For $0\le\gamma\le 1/2$ this does hold and we can proceed with the usual summation. By consolidating (12) into a complex exponential form, we have
\begin{equation}
\sum _{n=-\infty}^\infty\frac{(-1)^n}{n+\mu}e^{in\theta}=\frac{\pi}{\sin{\pi\mu}}e^{-i\mu\theta}
\end{equation}
for $\theta\in (-\pi ,\pi )$. Plugging this identity into the double integral representation of the Lerche-Newberger sum and separating the factors by integration variable gives us the result in (16). For $1/2<\gamma\le 1$, the quantity $\gamma (\phi -\theta )$ will extend beyond $(-\pi ,\pi )$. Indeed, if $\theta \in (\pi ,3\pi )$, then $\theta -2\pi\in (-\pi ,\pi )$ such that
\begin{equation}
\sum _{n=-\infty}^\infty\frac{(-1)^n}{n+\mu}e^{in\theta}=\frac{\pi}{\sin{\pi\mu}}e^{-i\mu (\theta -2\pi )}
\end{equation}
A similar identity holds for $\theta\in (-3\pi ,-\pi )$. If we restrict $\alpha ,\beta\in\Z$, then we can write the double integral representation as
\begin{equation}
\sum _{n=-\infty}^\infty\frac{(-1)^nA_{\alpha -\gamma n}({\bf x},{\bf y})A_{\beta +\gamma n}({\bf x},{\bf y})}{n+\mu}=\int _{-\pi}^\pi\int_{-\pi}^\pi f_\mu (\gamma (\phi -\theta ))p_\alpha (\theta )p_\beta (\phi )d\phi d\theta
\end{equation}
where $f _\mu(\theta)$ is the summation in (19), and $p_\alpha(\theta )$ is the $2\pi$-periodic integrand of $A_\alpha ({\bf x},{\bf y})$. In Figure 1, we decompose the region of integration according to the period that $\gamma (\phi -\theta )$ lies in. In the triangle region bounded by $\theta =-\pi$, $\phi =\pi$, and $\gamma (\phi -\theta )=\pi$, we have $\gamma (\phi -\theta )\ge \pi$ so the summation formula from (19) incurs an extra factor of $e^{2\pi\mu i}$ as shown in (20). Let this region be called $\Delta _+(\gamma )$. In the negation of that region, the triangle bounded by $\theta =\pi$, $\phi =-\pi$ and $\gamma (\phi -\theta )=-\pi$ which we call $\Delta _-(\gamma )$, the summation procedure incurs an extra factor of $e^{-2\pi\mu i}$. In this way, we can replace $f_\mu (\gamma (\phi -\theta ))$ by a complex exponential and consolidate the expression into an Anger function product and two integrals over $\Delta _\pm (\gamma )$:
\begin{align}
\sum _{n=-\infty}^\infty\frac{(-1)^nA_{\alpha -\gamma n}({\bf x},{\bf y})A_{\beta +\gamma n}({\bf x},{\bf y})}{n+\mu} &= \frac{\pi}{\sin{\mu\pi}}A_{\alpha +\gamma\mu}({\bf x},{\bf y})A_{\beta -\gamma\mu}({\bf x},{\bf y})+2\pi i e^{2i\pi\mu}\iint _{\Delta _+}p_{\alpha +\gamma\mu}(\theta )p_{\beta -\gamma\mu}(\phi ) d\theta d\phi \nonumber \\
	&- 2\pi i e^{-2i\pi\mu}\iint _{\Delta _-}p_{\alpha +\gamma\mu}(\theta )p_{\beta -\gamma\mu}(\phi ) d\theta d\phi
\end{align}

\begin{figure}
%\begin{center}\includegraphics[scale=0.5]{figure_1.png}\end{center}
\begin{center}\includegraphics[width=0.66\textwidth]{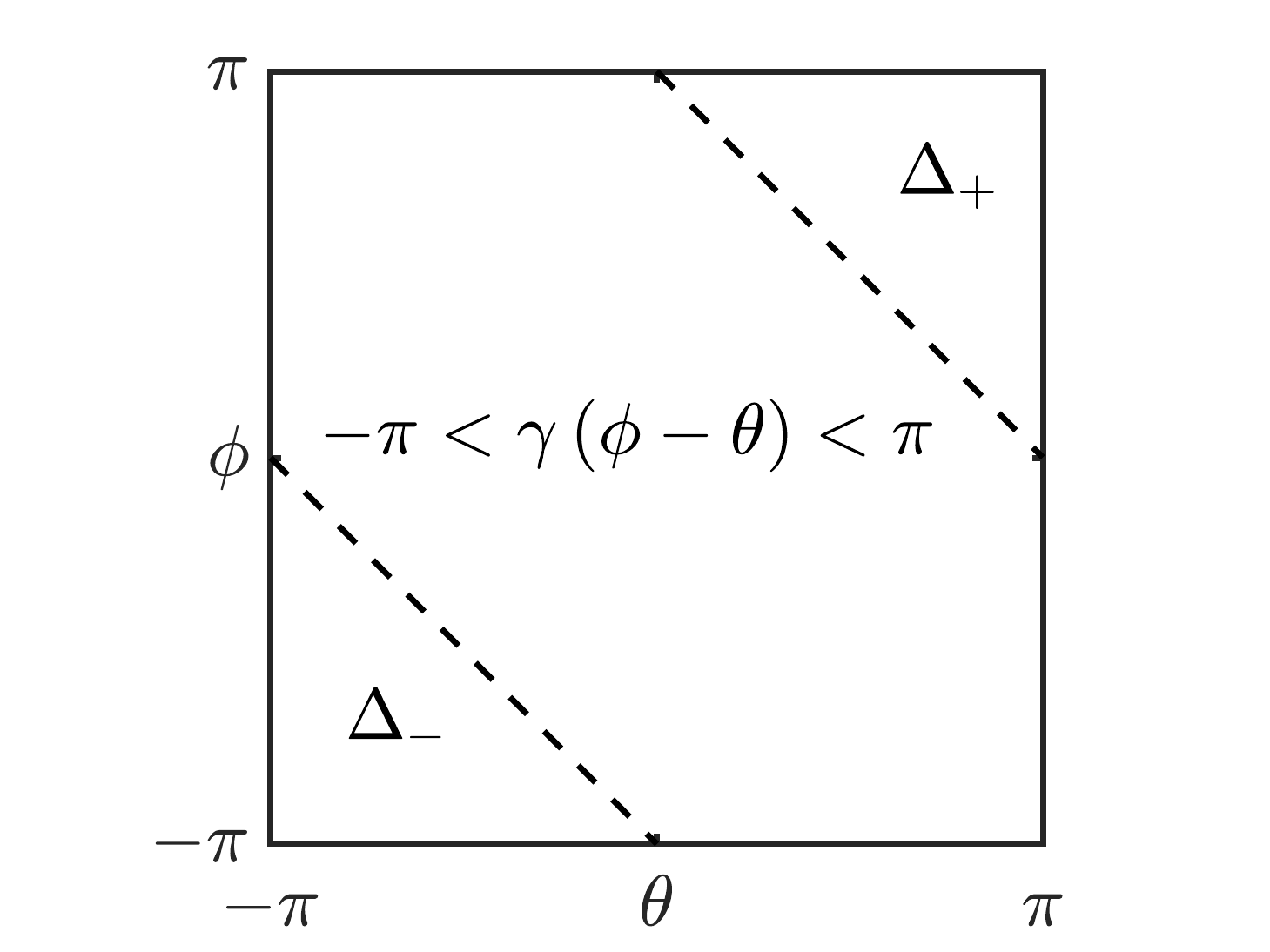}\end{center}
\caption{Division of region of integration in (22).}
\end{figure}

\section{Applications}

Lerche-Newberger type summations appear in many different areas of physics including plasma physics \cite{lerche74} and periodically driven quantum mechanical systems \cite{silveri17}. In particular, the following expression appears in Berns et. al. \cite{berns06} to describe the transition rate of a persistent qubit:
\begin{equation}
W(\epsilon ,A)=\frac{\Delta ^2}{2}\sum _n\frac{\Gamma _2J_n^2(x)}{(\epsilon -\omega n)^2+\Gamma _2^2}
\end{equation}
where $x=A/\omega$. The authors give an asymptotic treatment of this function, but a closed form expression is possible using (1). Using a partial fractions expansion, we can break the denominator into its linear factors. Let $\mu _\pm =(\epsilon\pm\Gamma i)/\omega$ such that
\begin{equation}
W(\epsilon ,A)=\frac{i\omega\Delta ^2}{4}\left[\sum _n\frac{J_n^2(x)}{n-\mu _+}-\sum _n\frac{J_n^2(x)}{n-\mu _-}\right]
\end{equation}
By letting $\alpha =\beta =0$ and $\gamma =1$, we can apply (1) to (24) to arrive at a closed form expression:
\begin{equation}
W(\epsilon ,A)=\frac{i\pi\omega\Delta ^2}{4}\left[\frac{J_{\mu _+}(x)J_{-\mu _+}(x)}{\sin\mu_+\pi}-\frac{J_{\mu _-}(x)J_{-\mu _-}(x)}{\sin\mu_-\pi}\right]
\end{equation}
We display plots of the expression in (25) for several selections of parameters in Figure 2. Now asymptotics need only be developed for the constituent Bessel functions rather than for the complicated sum in (23).

\begin{figure}
\begin{center}\includegraphics[width=0.75\textwidth]{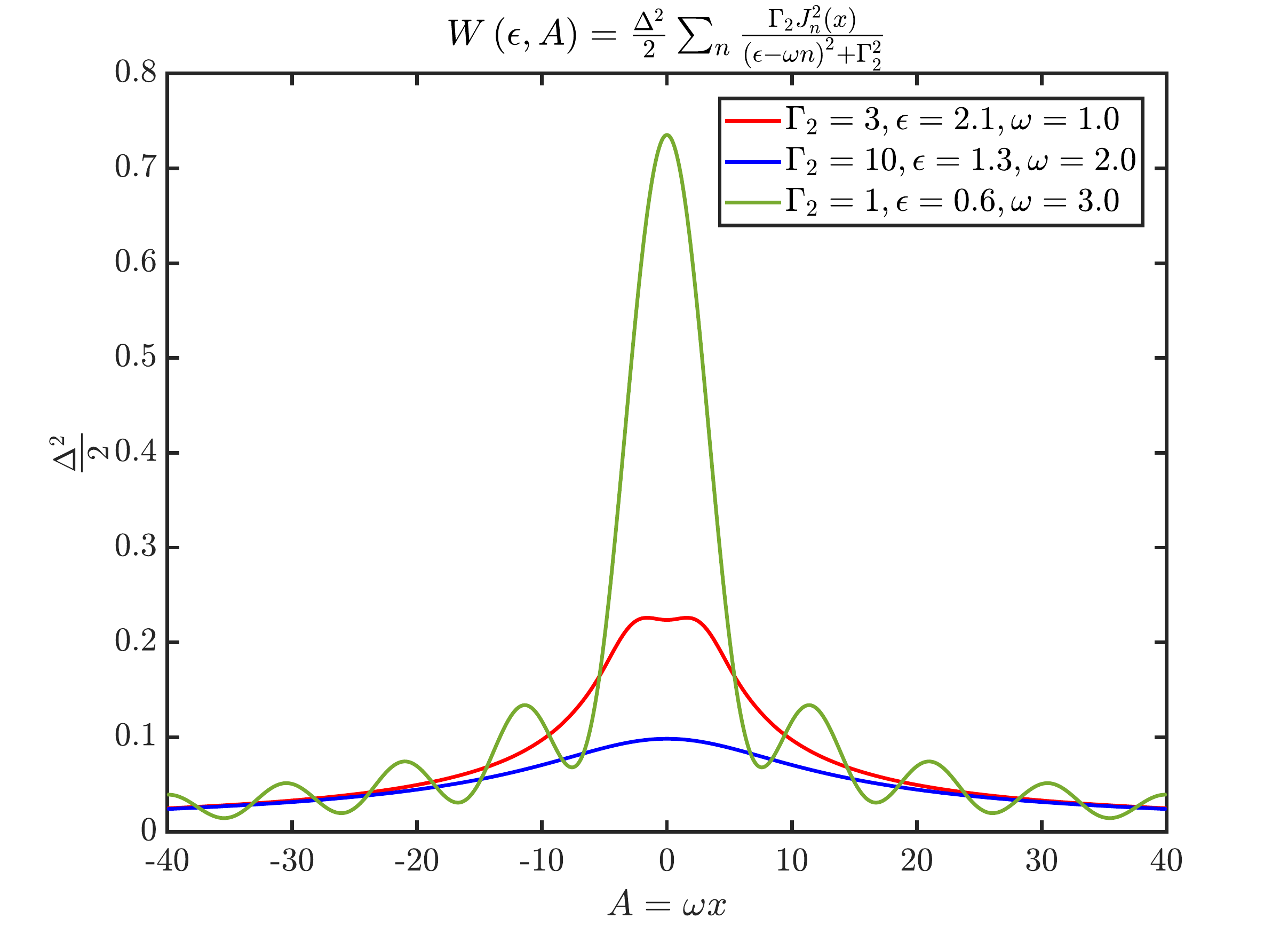}\end{center}
\caption{Plots of $W(\epsilon ,A)$ for a variety of parameter sets (\emph{Blue}) $\Gamma _2=3$, $\epsilon =2.1$, $\omega =0.07$  (\emph{Red}) $\Gamma _2=10$, $\epsilon =1.3$, $\omega =2.0$  (\emph{Green}) $\Gamma _2=1$, $\epsilon =0.6$, $\omega =3.0$.}
\end{figure}

If this system was modified by letting the qubit be modulated by several frequencies possibly out of phase where instead of the energy detuning being of the form $h(t)=\epsilon +\delta\epsilon (t)+A\cos{2\pi vt}$, we would include additional terms such that $h(t)=\epsilon +\delta\epsilon (t)+\sum _{k=1}^m\left( A_k\cos{2\pi vkt}+B_k\cos{2\pi vkt}\right)$. Then the Bessel functions in (23) would be replaced with generalized Bessel functions of the form $J_n({\bf x},{\bf y})$ where ${\bf x}=(A_1/\omega ,...,A_m/\omega )$ and ${\bf y}=(B_1/\omega ,...,B_m/\omega )$. After using another partial fractions decomposition, a variant of the generalized Lerche-Newberger summation would emerge:
\begin{equation}
\sum _{n=-\infty}^\infty\frac{J_n({\bf x},{\bf y})^2}{n+\mu}
\end{equation}
Since for generalized Bessel functions we have $J_{-n}({\bf x},{\bf y})\ne (-1)^nJ_n({\bf x},{\bf y})$, we must use a variant of (19). By representing $(-1)^n=e^{i\pi n}$ in (19), we can shift the piecewise continuous regions from $((2k-1)\pi ,(2k+1)\pi)$ to $(2k\pi ,2(k+1)\pi)$. Indeed, if $\theta\in (0,2\pi )$ we have
\begin{equation}
\sum _{n=-\infty}^\infty\frac{e^{in\theta}}{n+\mu}=\frac{\pi}{\sin\pi\mu}e^{i\mu\pi}e^{-i\mu\theta}
\end{equation}
If $\theta\in (-2\pi ,0)$, then a similar argument in (27) would show that
\begin{equation}
\sum _{n=-\infty}^\infty\frac{e^{in\theta}}{n+\mu}=\frac{\pi}{\sin\pi\mu}e^{-i\mu\pi}e^{-i\mu\theta}
\end{equation}
Using the typical integral interchange arguments, we can express (26) as the double integral:
\begin{equation}
\sum _{n=-\infty}^\infty\frac{J_n({\bf x},{\bf y})^2}{n+\mu}=\frac{1}{4\pi ^2}\int _{-\pi}^\pi\int _{-\pi}^\pi g_\mu (\phi +\theta )p_0(\theta )p_0(\phi )d\theta d\phi
\end{equation}
where $g_\mu (\theta )$ is the summation from (27) and (28) and $p_0(\theta )$ is the same as in the previous section. In Figure 3, we see how the discontinuity of $g_\mu (\theta )$ divides the region of integration into two subregions on either side of the curve $\phi +\theta =0$ which we refer to as $\tilde{\Delta}_\pm$. Therefore we can divide this into the following sum:
\begin{equation}
\sum _{n=-\infty}^\infty\frac{J_n({\bf x},{\bf y})^2}{n+\mu}=\frac{1}{4\pi\sin{\pi\mu}}\left[ e^{i\mu\pi}\iint _{\tilde{\Delta}_+}p_{-\mu}(\theta )p_{-\mu}(\phi )d\theta d\phi+e^{-i\mu\pi}\iint _{\tilde{\Delta}_-}p_{-\mu}(\theta )p_{-\mu}(\phi ) d\theta d\phi\right]
\end{equation}

\begin{figure}
%\begin{center}\includegraphics[scale=0.5]{figure_2.png}\end{center}
\begin{center}\includegraphics[width=0.66\textwidth]{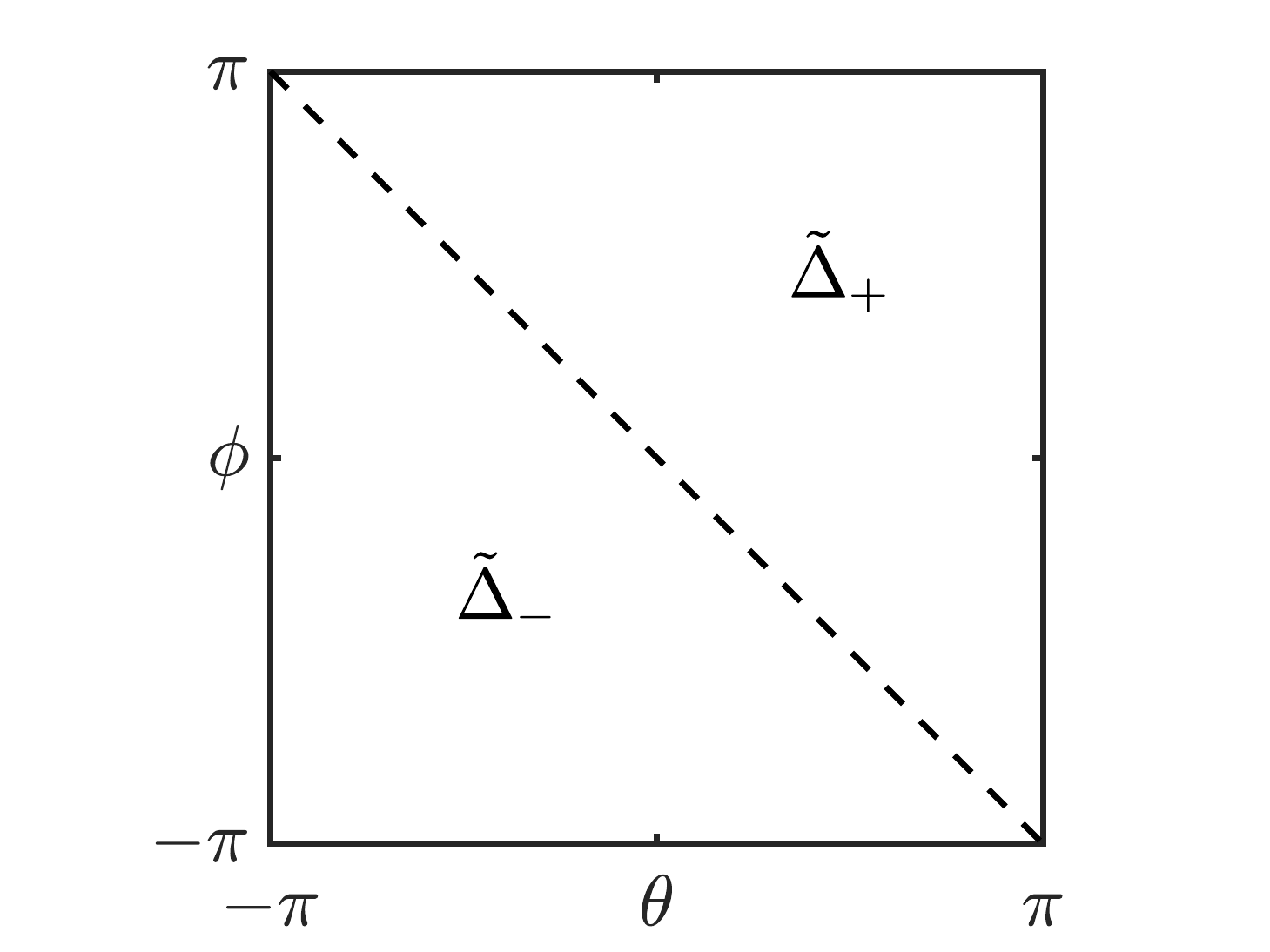}\end{center}
\caption{Division of region of integration in (30).}
\end{figure}

By splitting the complex exponentials into their trigonometric components, we find that this expression resolves based on the sum and difference of the two integral terms. Since $\tilde{\Delta}_+\cup\tilde{\Delta}_-=[-\pi,\pi]\times [-\pi ,\pi]$, it follows that
\begin{equation}\iint _{\tilde{\Delta}_+}p_{-\mu}(\theta )p_{-\mu}(\phi ) d\theta d\phi +\iint _{\tilde{\Delta}_-}p_{-\mu}(\theta )p_{-\mu}(\phi ) d\theta d\phi =4\pi ^2 A_{-\mu}({\bf x},{\bf y})^2
\end{equation}
Meanwhile, the difference of these two integrals can be combined into a single integral over the space by using the sign function $\text{sgn}(x)$ such that
\begin{equation}
\iint _{\tilde{\Delta}_+}p_{-\mu}(\theta )p_{-\mu}(\phi ) d\theta d\phi -\iint _{\tilde{\Delta}_-}p_{-\mu}(\theta )p_{-\mu}(\phi ) d\theta d\phi =\int _{-\pi}^\pi\int _{-\pi}^\pi p_{-\mu}(\theta )p_{-\mu}(\phi )\text{sgn}(\theta +\phi )d\theta d\phi
\end{equation}
and we represent this integral as $4\pi ^2B_{-\mu} ({\bf x},{\bf y})$ to compare with (31). Substituting (31) and (32) into (30) we arrive at the final expression:
\begin{equation}
\sum _{n=-\infty}^\infty\frac{J_n({\bf x},{\bf y})^2}{n+\mu}=2\pi\left(\cot (\pi\mu )A_{-\mu}({\bf x},{\bf y})^2+iB_{-\mu}({\bf x},{\bf y})\right)
\end{equation}
The asymptotics of $A_{-\mu}({\bf x},{\bf y})$ have been explored for ${\bf y}=0$ in several sources \cite{dattoli96} \cite{korsch06} and in more generality in Kuklinski and Hague \cite{kuklinski19}. The component $B_{-\mu}({\bf x},{\bf y})$ is a two-dimensional oscillatory integral over the region $[-\pi ,\pi ]\times [-\pi ,\pi ]$, however the standard results from Stein \cite{stein93} do not apply since the integrand has a discontinuity over the diagonal of this region contributed by $\text{sgn}(\phi +\theta )$ as seen in Figure 3. Regardless of the contributions from the discontinuity that may persist for all parameter choices $\mu ,{\bf x},{\bf y}$, we still see that $B_{-\mu}({\bf x},{\bf y})$ has the same bifurcation curves as $A_{-\mu}({\bf x},{\bf y})$ due to the similar structure of the integrand. We plot examples of these sums in Figure 4.

\begin{figure}
\begin{center}\includegraphics[scale=0.9]{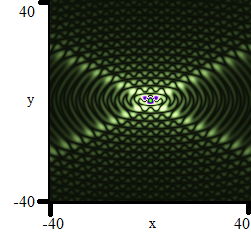}\includegraphics[scale=0.9]{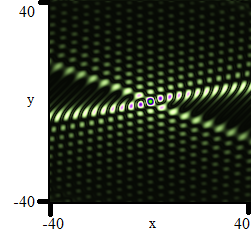}\includegraphics[scale=0.9]{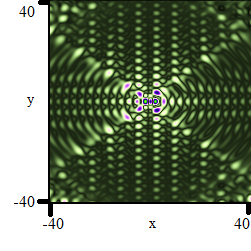}\end{center}

\begin{center}\includegraphics[scale=0.6]{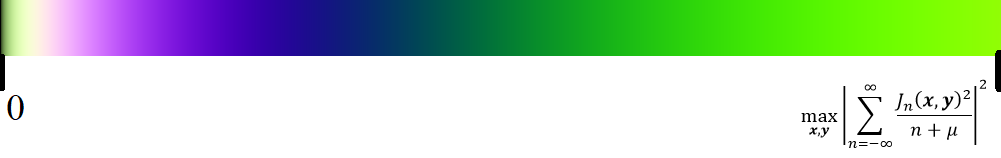}\end{center}
\caption{Plots of sum (33). (\emph{Left}) ${\bf x}=(x,y,0,...)$, ${\bf y}=0$, $\mu =0.5$ (\emph{Center}) ${\bf x}=(x,0,y,0,...)$, ${\bf y}=0$, $\mu =0.3$ (\emph{Right}) ${\bf x}=(0,x,y,0,...)$, ${\bf y}=0$, $\mu =1.7$.}
\end{figure}

\section{Conclusion}

In this document we extended the Lerche-Newberger formula to incorporate higher dimensional analogues to the Bessel functions. Since these generalized Bessel functions are only well defined for integer order, we instead opted to use generalize Anger functions in the summation. While the one-dimensional Lerche-Newberger formula generally preserves the Bessel function structure, the Anger function analogue does not behave as nicely since discontinuities in the equivalent integral expressions are introduced. Even for a relatively simple multi-dimensional as presented in the previous section, double integrals over discontinuous functions emerge that ostensibly cannot be reduced to well known functions using symmetry relations.

It is the hope of the authors that these results will contribute to future investigations of oscillatory quantum mechanical systems.

\section*{Data Availability Statement}

The data that support the findings of this study are available from the corresponding author upon reasonable request.

%\bibliography{biblio}

\begin{thebibliography}{10}
\providecommand{\url}[1]{#1}
\csname url@samestyle\endcsname
\providecommand{\newblock}{\relax}
\providecommand{\bibinfo}[2]{#2}
\providecommand{\BIBentrySTDinterwordspacing}{\spaceskip=0pt\relax}
\providecommand{\BIBentryALTinterwordstretchfactor}{4}
\providecommand{\BIBentryALTinterwordspacing}{\spaceskip=\fontdimen2\font plus
\BIBentryALTinterwordstretchfactor\fontdimen3\font minus
  \fontdimen4\font\relax}
\providecommand{\BIBforeignlanguage}[2]{{%
\expandafter\ifx\csname l@#1\endcsname\relax
\typeout{** WARNING: IEEEtran.bst: No hyphenation pattern has been}%
\typeout{** loaded for the language `#1'. Using the pattern for}%
\typeout{** the default language instead.}%
\else
\language=\csname l@#1\endcsname
\fi
#2}}
\providecommand{\BIBdecl}{\relax}
\BIBdecl

\bibitem{lerche08}
I.~Lerche and R.~Tautz, ``Kapteyn series arising in radiation problems,''
  \emph{Journal of physics A: Mathematical and theoretical}, vol.~41, no.~3, p.
  035202, 2008.

\bibitem{szapiel78}
S.~Szapiel, ``Mar{\'e}chal intensity criteria modified for circular apertures
  with nonuniform intensity transmission: Dini series approach,'' \emph{Optics
  letters}, vol.~2, no.~5, pp. 124--126, 1978.

\bibitem{watson22}
G.~Watson, \emph{A treatise on the theory of Bessel functions}.\hskip 1em plus
  0.5em minus 0.4em\relax Cambridge university press, 1922.

\bibitem{abramowitz64}
M.~Abramowitz and I.~A. Stegun, \emph{Handbook of mathematical functions with
  formulas, graphs, and mathematical tables}.\hskip 1em plus 0.5em minus
  0.4em\relax US Government printing office, 1964, vol.~55.

\bibitem{lerche07}
I.~Lerche and R.~C. Tautz, ``A note on summation of {K}apteyn series in
  astrophysical problems,'' \emph{The Astrophysical Journal}, vol. 665, no.~2,
  p. 1288, 2007.

\bibitem{kravchenko17}
V.~V. Kravchenko and S.~M. Torba, ``Asymptotics with respect to the spectral
  parameter and neumann series of bessel functions for solutions of the
  one-dimensional schr{\"o}dinger equation,'' \emph{Journal of Mathematical
  Physics}, vol.~58, no.~12, p. 122107, 2017.

\bibitem{linton06}
C.~Linton, ``Schl{\"o}milch series that arise in diffraction theory and their
  efficient computation,'' \emph{Journal of Physics A: Mathematical and
  General}, vol.~39, no.~13, p. 3325, 2006.

\bibitem{newberger82}
B.~S. Newberger, ``New sum rule for products of bessel functions with
  application to plasma physics,'' \emph{Journal of Mathematical Physics},
  vol.~23, no.~7, pp. 1278--1281, 1982.

\bibitem{buneman59}
O.~Buneman, ``Dissipation of currents in ionized media,'' \emph{Physical
  Review}, vol. 115, no.~3, p. 503, 1959.

\bibitem{lerche08a}
I.~Lerche, R.~Schlickeiser, and R.~Tautz, ``Comment on “a new derivation of
  the plasma susceptibility tensor for a hot magnetized plasma without infinite
  sums of products of bessel functions”[phys. plasmas 14, 092103 (2007)],''
  \emph{Physics of Plasmas}, vol.~15, no.~2, p. 092103, 2008.

\bibitem{silveri17}
M.~Silveri, J.~Tuorila, E.~Thuneberg, and G.~Paraoanu, ``Quantum systems under
  frequency modulation,'' \emph{Reports on Progress in Physics}, vol.~80,
  no.~5, p. 056002, 2017.

\bibitem{berns06}
D.~Berns, W.~Oliver, S.~Valenzuela, A.~Shytov, K.~Berggren, L.~Levitov, and
  T.~Orlando, ``Coherent quasiclassical dynamics of a persistent current
  qubit,'' \emph{Physical review letters}, vol.~97, no.~15, p. 150502, 2006.

\bibitem{kuklinski19}
P.~Kuklinski and D.~A. Hague, ``Identities and properties of multi-dimensional
  generalized bessel functions,'' \emph{arXiv preprint arXiv:1908.11683}, 2019.

\bibitem{dattoli96}
G.~Dattoli and A.~Torre, \emph{Theory and applications of generalized Bessel
  functions}.\hskip 1em plus 0.5em minus 0.4em\relax Arcane, 1996.

\bibitem{Hague_AES}
D.~A. {Hague}, ``Adaptive transmit waveform design using multitone sinusoidal
  frequency modulation,'' \emph{IEEE Transactions on Aerospace and Electronic
  Systems}, vol.~57, no.~2, pp. 1274--1287, 2021.

\bibitem{reiss62}
\BIBentryALTinterwordspacing
H.~R. Reiss, ``Absorption of light by light,'' \emph{Journal of Mathematical
  Physics}, vol.~3, no.~1, pp. 59--67, 1962. [Online]. Available:
  \url{https://doi.org/10.1063/1.1703787}
\BIBentrySTDinterwordspacing

\bibitem{paciorek94}
W.~Paciorek and G.~Chapuis, ``Generalized bessel functions in incommensurate
  structure analysis,'' \emph{Acta Crystallographica Section A: Foundations of
  Crystallography}, vol.~50, no.~2, pp. 194--203, 1994.

\bibitem{korsch06}
H.~Korsch, A.~Klumpp, and D.~Witthaut, ``On two-dimensional bessel functions,''
  \emph{Journal of Physics A: Mathematical and General}, vol.~39, no.~48, 2006.

\bibitem{bakker84}
M.~Bakker and N.~M. Temme, ``Sum rule for products of bessel functions:
  Comments on a paper by newberger,'' \emph{Journal of mathematical physics},
  vol.~25, no.~5, pp. 1266--1267, 1984.

\bibitem{lerche74}
I.~Lerche, ``A note on summing series of bessel functions occurring in certain
  plasma astrophysical situations,'' \emph{The Astrophysical Journal}, vol.
  190, pp. 165--166, 1974.

\bibitem{royden88}
H.~L. Royden and P.~Fitzpatrick, \emph{Real analysis}.\hskip 1em plus 0.5em
  minus 0.4em\relax Macmillan New York, 1988, vol.~32.

\bibitem{apostol51}
T.~M. Apostol \emph{et~al.}, ``On the lerch zeta function.'' \emph{Pacific
  Journal of Mathematics}, vol.~1, no.~2, pp. 161--167, 1951.

\bibitem{wilbraham48}
H.~Wilbraham, ``On a certain periodic function,'' \emph{Cambridge and Dublin
  Mathematical Journal}, vol.~3, pp. 198--201, 1848.

\bibitem{stein93}
E.~M. Stein, \emph{Harmonic Analysis: Real-Variable Methods, Orthogonality, and
  Oscillatory Integrals}.\hskip 1em plus 0.5em minus 0.4em\relax Princeton
  University Press, 1993.

\bibitem{hewitt79}
E.~Hewitt and R.~E. Hewitt, ``The gibbs-wilbraham phenomenon: an episode in
  fourier analysis,'' \emph{Archive for history of Exact Sciences}, vol.~21,
  no.~2, pp. 129--160, 1979.

\end{thebibliography}
%\bibliographystyle{IEEEtran}


\appendix
\section{Proof of Global Maximum Location}

In this section we prove that there exists a constant $M\in\N$ such that for all $n>M$, (7) has a global maximum at $\theta _n=\pi -\pi /(2n+1)$. We prove this for even $n$ so we replace $n$ with $2n$, but the result can easily be extended to $n$ odd. To do this we replicate an argument from Hewitt and Hewitt \cite{hewitt79}. Let $f(\theta )$ be the difference from (8) which we write here:
\begin{equation}
f(\theta )=\sum _{k=1}^{2n}\frac{(-1)^k}{k}\cos{k\theta}+\log\left(2\cos\frac{\theta}{2}\right)
\end{equation}
Using trigonometric identities, we can prove that the derivative of this function satisfies the following:
\begin{equation}
f'(\theta )=-\sin\left(\frac{4n+1}{2}\theta\right)/\left( 2\cos\left(\frac{\theta}{2}\right)\right)
\end{equation}
This gives us that the local extrema of $f$ satisfy $\theta _k=\frac{2\pi k}{4n+1}$, and by conducting a second derivative test we can conclude that $\theta _{2k}$ are the local maxima of $f$. We first will prove that $f(\theta _{2k+2})-f(\theta _{2k})>f(\theta _{2k})-f(\theta _{2k-2})$. To do this let us define a function $\omega$:
\begin{align}
\omega (t) &= \left[ f\left(\frac{4\pi (k+1/2)}{4n+1}+\frac{2t}{4n+1}\right) -f\left(\frac{4\pi (k-1/2)}{4n+1}+\frac{2t}{4n+1}\right)\right] \nonumber \\
	&- \left[ f\left(\frac{4\pi (k+1/2)}{4n+1}-\frac{2t}{4n+1}\right) -f\left(\frac{4\pi (k-1/2)}{4n+1}-\frac{2t}{4n+1}\right)\right]
\end{align}
Here, we restrict $1\le k\le n-1$. This function is of interest since if $\omega (\pi )>0$, then the inequality in question holds. Since $\omega (0)=0$, we need only show that $\omega '(t)>0$ for $0\le t\le \pi$. Indeed, $\omega '(t)$ takes the following form:
\begin{align}
\omega '(t) &= \frac{2\sin t}{4n+1}\left[\left(\frac{1}{\cos\left(\frac{2\pi (k+1/2)}{4n+1}+\frac{t}{4n+1}\right)}-\frac{1}{\cos\left(\frac{2\pi (k-1/2)}{4n+1}+\frac{t}{4n+1}\right)}\right)\right. \nonumber \\
	&- \left.\left(\frac{1}{\cos\left(\frac{2\pi (k+1/2)}{4n+1}-\frac{t}{4n+1}\right)}-\frac{1}{\cos\left(\frac{2\pi (k-1/2)}{4n+1}-\frac{t}{4n+1}\right)}\right)\right]
\end{align}
Using the identity $1/\cos(a+b)-1/\cos(a-b)=4(\sin{a}\sin{b})/(\cos{2a}+\cos{2b})$, we can rewrite the function as follows:
\begin{equation}
\omega '(t) = \frac{8\sin t\sin\left(\frac{\pi}{4n+1}\right)}{4n+1}\left[\frac{\sin\left(\frac{2\pi k+t}{4n+1}\right)}{\cos\left(\frac{4\pi k+2t}{4n+1}\right) +\cos\left(\frac{2\pi}{4n+1}\right)} -\frac{\sin\left(\frac{2\pi k-t}{4n+1}\right)}{\cos\left(\frac{4\pi k-2t}{4n+1}\right) +\cos\left(\frac{2\pi}{4n+1}\right)}\right]
\end{equation}
After combining these two fractions and conducting some elementary trigonometric manipulations, we can further simplify this derivative:
\begin{equation}
\omega '(t) = \frac{16\sin t\sin\left(\frac{\pi}{4n+1}\right)\sin\left(\frac{t}{4n+1}\right)\cos\left(\frac{2\pi k}{4n+1}\right)}{(4n+1)\left(\cos\left(\frac{4\pi k+2t}{4n+1}\right) +\cos\left(\frac{2\pi}{4n+1}\right)\right)}
\end{equation}
Due to the restrictions on $t$ and $k$, all of the terms in this fraction are positive and therefore the desired condition holds, namely that the difference between adjacent local maxima increases to the right.

We ultimately want to show that $\theta _{2n}$ is the global maxima. This can be accomplished by proving $f(\theta _2)>f(\theta _0)$ and using induction with the above result. Though this result appears true for all $n$, due to the unwieldy nature of $f(\theta)$ we instead opt to prove this asymptotically, that there exists some $M$ such that for all $n>M$ the inequality is satisfied. We do this by considering the asymptotic expansion of $f(\theta _2)-f(\theta _0)$:
\begin{equation}
f(\theta _2)-f(\theta _0)=-2\sum _{k=1}^{2n}\frac{(-1)^k}{k}\sin^2\left(\frac{2\pi k}{4n+1}\right) +\log\left(\cos\frac{2\pi}{4n+1}\right)
\end{equation}
If we can show that the leading term of the asymptotic expansion of (40) in $n$ is positive, then the proof follows. The logarithm term can be easily expanded using typical Taylor expansion arguments:
\begin{equation}
\log\left(\cos\frac{2\pi}{4n+1}\right) =\frac{1}{n^2}\left(-\frac{\pi ^2}{8}\right)+\frac{1}{n^3}\left(\frac{\pi ^2}{16}\right) +O(n^{-4})
\end{equation}
To resolve the summation term, we manipulate it into an endpoint Riemann sum. Recall that for a smooth function $f(x)$, the following asymptotic form holds \cite{wals37}:
\begin{equation}
\frac{1}{n}\sum _{k=1}^nf\left(\frac{k}{n}\right) =\int _0^1f(x)dx+\frac{1}{2n}\int _0^1f'(x)dx+\frac{1}{12n^2}\int _0^1f''(x)dx +O(n^{-4})
\end{equation}
We split the summation into positive and negative terms according to the $(-1)^k$ term:
\begin{equation}
\sum _{k=1}^{2n}\frac{(-1)^k}{k}\sin^2\left(\frac{2\pi k}{4n+1}\right) =\sum _{k=1}^n\frac{1}{2k}\sin ^2\left(\frac{4\pi k}{4n+1}\right) -\sum _{k=1}^n\frac{1}{2k-1}\sin ^2\left(\frac{2\pi (2k-1)}{4n+1}\right)
\end{equation}
We elaborate the asymptotic expansion of the even sum on the right hand side of (43); the odd sum is similar. By letting $f(x)=\sin ^2(\pi x)/x$, we rewrite this sum as
\begin{equation}
\sum _{k=1}^n\frac{1}{2k}\sin ^2\left(\frac{4\pi k}{4n+1}\right) =\frac{2}{4n+1}\sum _{k=1}^nf\left(\frac{4k}{4n+1}\right)
\end{equation}
The remainder of this proof is an arduous calculation of these sum up to order for large $n$ which we outline here. First, we expand the argument of $f$, $4k/(4n+1)$, being careful to recognize that $k$ is at the same order as $n$ so that $k/n=O(1)$. Next, we conduct a Taylor expansion of $f$ at $k/n$ throwing out higher order terms such that we are left with terms that look like $(k/n)^jf^{(j)}(k/n)$. By considering these the summands of a Riemann endpoint sum of function $g(x)=x^jf^{(j)}$, we use (42) to further expand these sums into a collection of integrals. These integrals have closed form expressions, and we can finish the problem off by multiplying the expansion of the sum on the right hand side of (44) by the expansion of the factor $2/(4n+1)$. This gives us the asymptotic expansions:
\begin{equation}
\sum _{k=1}^n\frac{1}{2k}\sin ^2\left(\frac{4\pi k}{4n+1}\right) =\frac{\pi x}{2}-\frac{1}{n^2}\left(\frac{\pi ^2}{24}\right) +\frac{1}{n^3}\left(\frac{5\pi ^2}{384}\right) +O(n^{-4})
\end{equation}
\begin{equation}
\sum _{k=1}^n\frac{1}{2k-1}\sin ^2\left(\frac{2\pi (2k-1)}{4n+1}\right) =\frac{\pi x}{2}+\frac{1}{n^2}\left(\frac{\pi ^2}{48}\right) -\frac{1}{n^3}\left(\frac{\pi ^2}{384}\right) +O(n^{-4})
\end{equation}
Here, $x=(-\text{Ci}(2\pi )+\gamma +\log (2\pi ))/2\pi$. Combining the expansions from (45) and (46) with the expansion of the logarithm in (41), we arrive at the final result, namely
\begin{equation}
f(\theta _2)-f(\theta _0)=\frac{\pi ^2}{32n^3}+O(n^{-4})
\end{equation}
and since the leading quantity is positive, the rest of the proof follows.
 
\end{document}